\documentclass{amsart}

\usepackage{amsmath,amsthm,amssymb}

\newtheorem{theorem}{Theorem}[section]
\newtheorem{lemma}[theorem]{Lemma}

\newtheoremstyle{spec}
  {9pt}
  {9pt}
  {\itshape}
  {}
  {\bfseries}
  {.\,\quad }
  {2pt}
  {}
\theoremstyle{spec}
\newtheorem{utheorem}{Theorem}

\newtheoremstyle{pspec}
  {9pt}
  {9pt}
  {\itshape}
  {}
  {\bfseries}
  {$\varepsilon$.\,\: }
  {2pt}
  {}
\theoremstyle{pspec}
\newtheorem{ptheorem}{Theorem}

\newtheoremstyle{cspec}
  {9pt}
  {9pt}
  {\itshape}
  {}
  {\bfseries}
  {$\sigma$.\, }
  {2pt}
  {}
\theoremstyle{cspec}
\newtheorem{ctheorem}{Theorem}

\theoremstyle{definition}
\newtheorem{definition}[theorem]{Definition}

\theoremstyle{remark}

\numberwithin{equation}{section}

\newcommand{\n}[1]{\left\|#1\right\|}
\newcommand{\m}[1]{\left|#1\right|}

\renewcommand{\epsilon}{\varepsilon}

\begin{document}

\title[Eigenvalues, Pseudospectra and Conditionspectra]{Comparative results on Eigenvalues, Pseudospectra and Conditionspectra}
\author[suku]{D. Sukumar}



\date{\today}



\begin{abstract}
Conditionspectrum measures the computational stability of solving a linear system. In this paper, ten theorems involving $\epsilon$-conditionspectrum are presented. All these theorems generalize a well known eigenvalue theorem and simultaneously compare with an appropriate pseudospectra theorem. Our organizing principle is that each conditionspectrum result precisely reduces to the corresponding eigenvalue theorem when $\epsilon=0$. The format of each conditionspectrum result is similar to the pseudospectrum result for easy comparison.
Each condition spectrum is formatted similar to pseudospectrum result for the easy comparison.
\end{abstract}

\maketitle


\bibliographystyle{amsplain}



The idea of \textit{condition spectrum} was conceived while trying to realize pseudospectra as a special case of generalized spectra defined by Ransford  in \cite{Rans}. The axiomatic approach of his study compelled S.H.Kulkarni and the author in \cite{suku2}, to modify slightly the underlying basic set so as to satisfy the axioms of Ransford. The newly defined condition spectrum in itself has interesting properties which were published in \cite{suku2,suku3}. In this article, its connection with the usual spectrum (eigenvalue) and pseudospectrum are discussed. Only those results on condition spectrum that have exemplary analogy with the other two spectra are chosen for presentation. 

Each theorem is arranged in triplet, labeled with the same number followed by relevant suffix. The first part is a result about the usual spectrum. The second part is a generalization of the first result using pseudospectrum, while the third part, the main contribution of this article, generalizes the first part result in terms of condition spectrum. For example, Theorem $1$, $1\varepsilon$ and $1\sigma$, denotes the first theorem for usual spectrum, pseudospectrum and condition spectrum respectively. This arrangement makes the reader to compare and comprehend the results conveniently. Note that all the $\varepsilon$ and $\sigma$ theorems reduces to corresponding usual spectrum (eigenvalue) theorem, as $\varepsilon$  tends to zero. This arrangement of results are essentially following the article \cite{Mark}.

\section{Introduction}
Let $A$ be a $N\times N$ matrix with complex entries., ie $A \in M_N(\mathbb C)$. Let $\lambda$ and $z$	denote complex numbers. 
\begin{definition}[Condition spectrum] Let $0<\epsilon<1$. The $\epsilon$-condition spectrum of a matrix $A$ is defined as  
$$\sigma_\epsilon(A):=\left\{z\in \mathbb C:z-A \text{ is not invertible or} \n{(z-A)^{-1}}\n{z-A}\geq
\frac{1}{\epsilon}\right\}.$$ 
\end{definition}
Equivalently we can write it as $\sigma_\epsilon(A):=\left\{z\in \mathbb C: \n{(z-A)^{-1}}\n{z-A}\geq
\frac{1}{\epsilon}\right\}$ with the convention that $\n{(z-A)^{-1}}\n{z-A}=\infty$ when $z-A$ is not invertible. The name of the spectrum and the range of $\epsilon$ are justified from the definition, as it related to the condition number of $A-z$.

For the sake of completeness we give the definition of pseudospectrum which is used to understand the behavior of non-normal matrices. The properties and applications are available in \cite{trefbook}.
\begin{definition}[Pseudospectrum] Let $\epsilon>0$. The $\epsilon$-pseudospectrum of a matrix $A$ is defined as  
$$\Lambda_\epsilon(A):=\left\{z\in \mathbb C:z-A \text{ is not invertible or} \n{(z-A)^{-1}}\geq
\epsilon\right\}.$$ 
\end{definition}
Consider solving the system of equations $Ax-\lambda x=b$. The eigenvalues of $A$ says the points in the complex plane at which $A-\lambda$ is not invertible and hence deals about the uniqueness of the solution. The pseudospectra of $A$ conveys the points at which $\n{(A-\lambda)^{-1}}$ becomes very large and hence deals with the computational aspects of the solution. The condition spectra captures the points in the complex plane at which the condition number $\n{A-\lambda}\n{(A-\lambda)^{-1}}$ becomes huge and hence elaborates the computational stability aspect of deriving the solution.

It is proved in \cite{suku2} that the condition spectrum is a non empty, compact, perfect set (no isolated points) and always contains the usual spectrum.  In \cite{suku3}, a sufficient condition for a function to be almost multiplicative is given using conditionspectrum. Here we list two results which will be referred often subsequently. The interested reader can refer \cite{suku2} for the proofs.
\begin{lemma}
For every $0<\epsilon<1$, $\sigma(A)\subseteq \sigma_\epsilon(A)$.
\end{lemma}

\begin{lemma}
For every $0<\epsilon<1$, $\sigma_{\epsilon}(A)$ is compact.
\end{lemma}

The above results are proved in a general setting, when $A$ is a Banach algebra element. On the contrast, the proof of the following equivalent conditions mainly relies on finite dimensionality. 
\begin{theorem}
The following definitions are equivalent.
\begin{enumerate}
\item $\sigma_\epsilon(A)=\{z\in \mathbb C:\n{(z-A)^{-1}}\n{z-A}\geq
\epsilon^{-1}\}$
\item $B=\{z\in \mathbb C:\exists\,
u\in {\mathbb C}^n\,\text{with}\, \n{u}=1 \,\text{s.t}\, \n{(z-A)u}\leq 
\epsilon\n{(z-A)}\}$
\item $C=\{z\in \mathbb C:z\in\sigma(A+E) \,\text{for
some}\, E \,
\text{with}\n{E}\leq\epsilon\n{(z-A)}\}$
\end{enumerate}
\end{theorem}

%

\begin{proof}
First observe that eigenvalues of $A$ are there in all the three sets. For that, 
suppose $z$ is an eigenvalue then
\begin{itemize}
\item $z\in \sigma_\epsilon(A)$ by the convention $\n{z-A}\n{(z-A)^{-1}}=\infty$.
\item $z\in B$ by taking $u$ as any normalized eigenvector
corresponding to the eigenvalue $z$.
\item $z\in C$ by taking $E=0$.
\end{itemize}
As eigenvalues belongs to all sets, it is enough to prove, the elements of one set (excluding eigenvalues) belongs to other sets.
The following three implications together establish the proof of the theorem.

$(1\Rightarrow 2):$
Suppose $z\in \sigma_\epsilon(A)$ and $z\notin \sigma(A)$ then $\n{z-A}\n{(z-A)^{-1}}\geq \epsilon^{-1}$. Since unit sphere
in finite dimensional space is compact, there exists an element $u\in\mathbb C^n$ such that
$\n{u}=1$ and $$\n{(z-A)^{-1}u}=\n{(z-A)^{-1}}$$  Define $$\widetilde{u}=(z-A)^{-1}u,$$ so
that $\n{\widetilde{u}}=\n{(z-A)^{-1}}$ and $(z-A)\widetilde{u}=u$.
Let $\widehat{u}=\dfrac{\widetilde{u}}{\n{\widetilde{u}}}$, then
$$\n{(z-A)\widehat{u}}=
\dfrac{\n{(z-A)\widetilde{u}}}{\n{\widetilde{u}}}=\dfrac{\n{u}}{\n{\widetilde{u}}}=
\dfrac{1}{\n{(z-A)^{-1}}}\leq
\epsilon\n{z-A}$$ and hence $z\in B$.

$(2\Rightarrow 3):$
Suppose $z\in B$ then there exists $u\in \mathbb C^n$ with $\n{u}=1$ such that
$$\n{(z-A){u}}\leq \epsilon\n{(z-A)}.$$ Let $v\in \mathbb C^n$ be a unit vector
satisfying $(A-z)u=\hat{\epsilon}v$ with $\hat{\epsilon}\leq \epsilon\n{z-A}$.
Let $w$ be a vector with $\n{w}=1$ such that $w^*u=1$. With such a $w$ we can
write $$zu=Au-\hat{\epsilon}vw^*u=(A-\hat{\epsilon}vw^*)u$$which means that
$z\in \sigma(A+E)$ for $E=-\hat{\epsilon}vw^*$ satisfying $\n{E}\leq
\hat{\epsilon}$. 

$(3\Rightarrow 1):$
Suppose $z\in\sigma(A+E)$ for some $E$
with$\n{E}\leq \epsilon\n{(z-A)}$. Then there exists a unit vector $v\in \mathbb
C$ such that $(A+E)v=zv$. 
By rearranging and
inverting
\begin{align*}
zv&=Av+Ev\\
(z-A)v&=Ev\\
v&=(z-A)^{-1}Ev
\end{align*}
and thus we have
$$1=\n{v}=\n{(z-A)^{-1}Ev}\leq\n{(z-A)^{-1}}\n{E}\leq\n{(z-A)^{-1}}\epsilon\n{(z
-A)},$$ implies that $\n{(z-A)^{-1}}\n{(z-A)}\geq{\epsilon}^{-1}$ and hence $z\in \sigma_\epsilon(A)$.\qedhere 
\end{proof}
\section{Results}
Ten results that are generalizing the theorems of usual spectrum are given. 
Only proofs of conditionspectrum results, that is theorems with the suffix $\sigma$, are given. Proofs of the usual spectrum are available in any standard book on linear algebra and the proofs of pseudospectra results are available in \cite{Mark}. Note that each condition spectrum result tends to the corresponding eigenvalue theorem in the limiting case as $\epsilon$ tends to 0.
\begin{utheorem}
$A$ is singular $\Longleftrightarrow$ $0\in\sigma(A)$.
\end{utheorem}
\begin{ptheorem}
$\n{A^{-1}}\geq \epsilon^{-1}$ $\Longleftrightarrow$ $0\in\Lambda_{\epsilon}(A)$ when $A$ is not singular. 
\end{ptheorem}
\begin{ctheorem}
$\n{A}\n{A^{-1}}\geq \epsilon^{-1}$ $\Longleftrightarrow$ $0\in\sigma_{\epsilon}(A)$ when $A$ is not singular. 
\end{ctheorem}
\begin{proof}
Follows from definition of condition spectrum.
\end{proof}
The bounds for the spectral radius of the condition spectrum.	
 \begin{utheorem}
 $\lambda\in\sigma(A)\Rightarrow \m{\lambda}\leq\n{A}$
 \end{utheorem}
 \begin{ptheorem}
 $\lambda\in\Lambda_{\epsilon}(A)\Rightarrow \m{\lambda}\leq\n{A}+\epsilon$
 \end{ptheorem}
 \begin{ctheorem}\label{2}
 $\lambda\in\sigma_{\epsilon}(A)\Rightarrow \m{\lambda}\leq\dfrac{1+\epsilon}{1-\epsilon}\n{A}$
\end{ctheorem}
 \begin{proof}
Let $\lambda\in \sigma_{\epsilon}(A)$. If $\m{\lambda}\leq\n{A}$, then clearly $\m{\lambda}\leq\dfrac{1+\epsilon}{1-\epsilon}\n{A}$. 

Suppose $\m{\lambda}>\n{A}$ then $\lambda - A$ is invertible and $\n{(\lambda - A)^{-1}}\le\dfrac{1}{\m{\lambda}-\n{A}}$.
Using this in the definition of condition spectrum we have $$1\le \epsilon \n{(\lambda - a)^{-1}}\n{\lambda - a}\leq
\epsilon\frac{\m{\lambda}+\n{A}}{\m{\lambda}-\n{A}}.$$ On simplification,
\begin{equation*}
\m{\lambda}\leq \dfrac{1+\epsilon}{1-\epsilon}\n{A}.\hfill \qedhere 
\end{equation*}
\end{proof}

\begin{utheorem}
$A$ has $N$ distinct eigenvalues $\Rightarrow$ $A$ is diagonalizable.
\end{utheorem}
\begin{ptheorem}
$\Lambda_\epsilon(A)$ has $N$ distinct components $\Rightarrow$ $A$ is diagonalizable.
\end{ptheorem}
\begin{ctheorem}
$\sigma_\epsilon(A)$ has $N$ distinct components $\Rightarrow$ $A$ is diagonalizable.
\end{ctheorem}
\begin{proof}
Refer 
\cite{suku2} for a proof.
\end{proof}

There are two components in condition number. One is norm of the matrix and the other is norm of its inverse. The following inequality estimates the norm of the inverse. 
\begin{lemma}\label{lemma1}
Let $\mathcal A$ be a Banach algebra. Let $a,b\in \mathcal A$. If $a$ is invertible in $\mathcal A$ and $b$ is not invertible, then 
\begin{equation}
\dfrac{1}{\n{a^{-1}}}\label{1}\leq \n{a-b}.
\end{equation}
\end{lemma}
The next result gives a relation between the norm of the resolvent with condition spectrum.
\begin{utheorem}
$\n{(z-A)^{-1}}\geq\dfrac{1}{d(z,\sigma(A))}$
\end{utheorem}
\begin{ptheorem}
$\n{(z-A)^{-1}}\geq\dfrac{1}{d(z,\Lambda_{\epsilon}(A))+\epsilon}$
\end{ptheorem}

\begin{ctheorem}
$\n{(z-A)^{-1}}\geq\dfrac{1}{d(z,\sigma_{\epsilon}(A))+\frac{2\epsilon}{1-\epsilon}\n{A}}$
\end{ctheorem}
\begin{proof}
If $z\in\sigma_{\epsilon}(A)$ then the inequality is immediate from the definition of condition spectrum. So assume 
$z\notin\sigma(A)$ that is $z-A$ is invertible. Since $\sigma_{\epsilon}(A)$ is compact, 
we can choose $\lambda\in \sigma_\epsilon(A)$ such that  $\m{\lambda-z}=d(z,\sigma_{\epsilon}(A))$. As $\lambda\in \sigma_\epsilon(A)$ there exist a matrix $E$ such that $\n{E}\leq \epsilon \n{\lambda-A}$ such that $\lambda\in\sigma(A+E)$. This implies $\lambda-A-E$ is not invertible and so $z-(A+E+z-\lambda)$ also. Now letting $a=z-A$ and $b=z-(A+E+z-\lambda)$ in Lemma \ref{lemma1}, we get
\begin{align*} 
\dfrac{1}{\n{(z-A)^{-1}}}&\leq \n{z-A-[z-(A+E+z-\lambda)]}=\n{E+z-\lambda}\\
&\leq\m{z-\lambda}+\n{E}\\&\leq d(z,\sigma_{\epsilon}(A))+\epsilon \n{\lambda-A}\leq d(z,\sigma_{\epsilon}(A))+\frac{2\epsilon}{1-\epsilon} \n{A}
\end{align*}
Hence it establishes the required inequality.
\end{proof}
The spectrum is invariant under similarity transformation. Similarly the condition spectrum is also preserved under certain similarity transformation. This result will be helpful in the computational designing of conditionspectra. We use the standard notation $\kappa(S)=\n{S}\n{S^{-1}}$ to denote the condition number of a matrix $S$.
\begin{utheorem}
$A=SBS^{-1}\Rightarrow \sigma(A)=\sigma(B)$
\end{utheorem}
\begin{ptheorem}
$A=SBS^{-1}\Rightarrow \Lambda_{\epsilon}(A)\subseteq\Lambda_{\kappa(S)\epsilon}(B)$
\end{ptheorem}
\begin{ctheorem}
$A=SBS^{-1}\Rightarrow \sigma_\epsilon(A)\subseteq\sigma_{\kappa(S)^{2}\,\epsilon}(B)$ whenever $\kappa(S)^2\epsilon <1$.
\end{ctheorem}
\begin{proof}Let  $z\in\sigma_{\epsilon}(A)$ then 
\begin{align*}
\frac{1}{\epsilon}\leq \n{z-A}\n{(z-A)^{-1}}
&=\n{zSS^{-1}-SBS^{-1}}\n{(zSS^{-1}-SBS^{-1})^{-1}}\\
&\leq\left(\n{S}\n{S^{-1}}\right)^{2}\n{z-B}\n{(z-B)^{-1}}\\
&\leq \kappa(S)^2\n{z-B}\n{(z-B)^{-1}} \qedhere
\end{align*}
\end{proof}
From this it is clear that similarity transformation through a matrix $S$ with condition number 1, that is $\kappa(S)=1$, preserves the condition spectrum, (since $\sigma_{\epsilon}(A)\subseteq \sigma_{\epsilon}(B)$ and $\sigma_{\epsilon}(B)\subseteq \sigma_{\epsilon}(A)$). In particular, similarity transformation through orthogonal and unitary matrices preserves the condition spectrum. 
%

%
The next two results calculate the transient behavior of $A$ from the knowledge of its condition spectrum. But they do not completely describe the nature, as pointed out in \cite{ransfordpseudo}. If one of the eigenvalue is bigger than 1 then the powers of $A$ blows up to infinity. Condition spectrum also behaves in the same way. 
\begin{utheorem}
\begin{displaymath}
\max_{\lambda\in\sigma(A)}\m{\lambda}>1 \Rightarrow \sup_{k\geq 0}\n{A^k}=\infty
\end{displaymath}
\end{utheorem}
\begin{ptheorem}
$$\max_{\lambda\in\Lambda_\epsilon(A)}\m{\lambda}>1+M
\epsilon \Rightarrow \sup_{k>0}\n{A^k} > M$$
\end{ptheorem}
\begin{ctheorem}
$$\max_{\lambda\in\sigma_\epsilon(A)}\m{\lambda}>\dfrac{1+M^2 \epsilon}{1-M
\epsilon} \Rightarrow \sup_{k\geq0}\n{A^k} > M\quad \text{whenever}\quad M\leq \frac{1}{\epsilon}$$
\end{ctheorem}
\begin{proof}
It is easy to establish the result for two simple and extreme cases. First one is the case when $M<1$. Here the result is immediate as $\n{A^0}=1$. The other case is when there is an eigenvalue $\lambda$ of $A$ with $\m{\lambda}>1$. Since $\lambda$ is an eigenvalue of $A$, $\lambda^k$ is an eigenvalue of $A^k$ for all $k\in \mathbb N$ and hence $\m{\lambda^k}\leq \n{A^k}$ for all $k$. This implies $\sup_{k\geq 0}\n{A^k}$ is infinity and hence the result. 

Now, excluding the above mentioned simple cases, we prove the theorem by negation. That is, assume $M\geq 1$ and no eigenvalue of $A$ has absolute value greater than 1. Suppose $\sup_{k\geq 0}\n{A^k}\leq M$,
we will prove $$\max_{\lambda\in \sigma_{\epsilon}(A)}\leq \frac{1+M^2\epsilon}{1-M\epsilon}.$$
%
Let $\lambda\in \sigma_\epsilon(A)\setminus\sigma(A)$. It is clear that if $\m{\lambda}\leq M$ then $\m{\lambda} \leq \frac{1+M^2\epsilon}{1-M\epsilon}$ whenever $M< \frac{1}{\epsilon}$. 
When $\m{\lambda}>M$, as $\lambda-A$ is invertible, we get 
\begin{align*}
\n{(\lambda-A)^{-1}}\leq \dfrac{1}{\m{\lambda}}\sum_{k=0}^{\infty}\dfrac{\n{A^{k}}}{\m{\lambda}^{k}}\leq \dfrac{M}{\m{\lambda}}\sum_{k=0}^{\infty}\dfrac{1}{\m{\lambda}^{k}}=\dfrac{M}{\m{\lambda}}\left(\dfrac{1}{1-\dfrac{1}{\m{\lambda}}}\right)=\dfrac{M}{\m{\lambda}-1}.  
\end{align*}
Combining the above inequality with $\lambda\in\sigma_\epsilon(A)$ gives $\dfrac{1}{\epsilon \n{\lambda-A}}\leq\dfrac{M}{\m{\lambda}-1}$. On simplification, $\m{\lambda}\leq \dfrac{1+M^2 \epsilon}{1-M
\epsilon}$ provided $M<\dfrac{1}{\epsilon}$. As $\lambda$ is arbitrary in $\sigma_\epsilon(A)\setminus \sigma (A)$ $$\max_{\lambda\in\sigma_\epsilon(A)}\m{\lambda}\leq\dfrac{1+M^2 \epsilon}{1-M
\epsilon},$$ and that proves the claim.
\end{proof}
\begin{utheorem}
$\lambda\in \sigma(A) \Rightarrow \n{A^k}\geq \m{\lambda}^k$ for all $k$
\end{utheorem}
\begin{ptheorem}
$\lambda\in \Lambda_\epsilon(A) \Rightarrow \n{A^k}\geq \m{\lambda}^k-\dfrac{k\epsilon\n{A}^{k-1}}{1-k\epsilon/\n{A}}$ for all $k$ such that $k\epsilon <\n{A}$
\end{ptheorem}
\begin{ctheorem}
$\lambda\in \sigma_\epsilon(A) \Rightarrow \n{A^k}\geq \m{\lambda}^k-\dfrac{ks\n{A}^{k-1}}{1-ks/\n{A}}$ for all $k$ such that $(2k+1)\epsilon <1$ where $s=\dfrac{2\epsilon}{1-\epsilon}\n{A}$
\end{ctheorem}
\begin{proof}
Let $\lambda \in \sigma_\epsilon(A)$. Pick $\n{E}\leq \epsilon \n{\lambda-A}$ such that $\lambda\in\sigma(A+E)$. Then $\m{\lambda^k}\leq\n{(A+E)^k}$ which implies 
\begin{align*}
\n{A}^k&\geq\m{\lambda}^k- \binom{k}{1} \n{A}^{k-1}\n{E}-\binom{k}{2} \n{A}^{k-2}\n{E}^{2}-\dots\\
&\geq \m{\lambda}^k-k \n{A}^{k-1}\epsilon \n{\lambda-A}\left(1+ \frac{k\epsilon \n{\lambda-A}}{\n{A}}+ \left(\frac{k\epsilon \n{\lambda-A}}{\n{A}}\right)^2+\dots\right). 
\end{align*}
Provided $k\epsilon \n{\lambda-A}\leq \n{A}$, the series in this last equation converges, giving
\[\n{A}^k\geq \m{\lambda}^k-\frac{k \n{A}^{k-1}\epsilon \n{\lambda-A}}{1-\frac{k\epsilon \n{\lambda-A}}{\n{A}}}\geq \m{\lambda}^k-\frac{ks \n{A}^{k-1}}{1-\frac{ks}{\n{A}}}\] 
with $s=\dfrac{2\epsilon}{1-\epsilon}\n{A}$. Also note that $(2k+1)\epsilon <1$ implies $k\epsilon \n{\lambda-A}\leq \n{A}$ by Theorem \ref{2}$\sigma$.
\end{proof}

The following result is analogous with Gerschgorin's theorem which locate the spectrum in the complex plane. In our case it locates the condition spectrum. We denote the closed ball with center $a$ and radius $r$ by $D(a,r)$. Let $d_j$ be the diagonal entry $a_{jj}$ of $A$ and $r_j=\sum_{i=1 \& i\neq j}^N\m{a_{ij}}$.
\begin{utheorem}
$\sigma(A)\subseteq  \bigcup_{j=1}^N D(d_j,r_j)$
\end{utheorem}
\begin{ptheorem}
$\Lambda_{\epsilon}(A)\subseteq  \bigcup_{j=1}^N D(d_j,r_j+\sqrt{N}\epsilon)$
\end{ptheorem}
\begin{ctheorem}
$\sigma_\epsilon(A)\subseteq  \bigcup_{j=1}^N D\left(d_j,r_j+\sqrt{N}\frac{2\epsilon}{1-\epsilon}\n{A}\right)$
\end{ctheorem}
\begin{proof}
By equivalent definition,  $\lambda\in\sigma_\epsilon(A)$ if and only if $\lambda\in\sigma(A+E)$ for some $E$ with $\n{E}\leq \epsilon\n{\lambda-A}$. Applying the Greschgorin's theorem to $A+E$ we get $$\lambda\in\bigcup_jD\left(d_j+e_{jj},r_j+\sum_{k\neq j}{\m{e_{jk}}}\right)$$ where $e_{ij}$ are the entries of $E$. The above balls are contained in balls centered at $d_j$ with radius $r_j+\sum_{k=1}^{N} \m{e_{jk}}$. Now note that, with $E_j$ as a matrix whose $j$th row is same as $j$th row of $E$ and remaining rows are zero, $$\sum_{k=1}^{N} \m{e_{jk}}\leq \n{E_j}_\infty\leq \sqrt(N)\n{E_j}_2\leq \sqrt{N}\n{E}_2.$$ We know that $\n{E}\leq \epsilon\n{\lambda-A}\leq \epsilon \left(\frac{1+\epsilon}{1-\epsilon}\n{A}+\n{A}\right)\leq \frac{2\epsilon}{1-\epsilon}\n{A}$. Hence using this information and the previous inequality we get $$\lambda\in\bigcup_jD\left(d_j,r_j+\sqrt{N}\frac{2\epsilon}{1-\epsilon}\n{A}\right)$$
Since $\lambda$ is arbitrary, this proves the required result.
\[\sigma_\epsilon(A)\subseteq  \bigcup_{j}D\left(d_j,r_j+\sqrt{N}\frac{2\epsilon}{1-\epsilon}\n{A}\right)\qedhere\]
\end{proof}

Numerical range is equally explored like spectrum because of its computability. The next result 
connects numerical range and condition spectrum. We write $conv(S)$ for the convex hull in $\mathbb C$ of a set $S\subseteq \mathbb C$. The notion $S\setminus \epsilon$-border means the set of points $z\in \mathbb C$ such that $D(z,\epsilon)\subseteq S$
\begin{utheorem}
$ W(A)\supseteq conv(\sigma(A))$
\end{utheorem}
\begin{ptheorem}
$ W(A)\supseteq conv(\Lambda_\epsilon(A))\setminus \epsilon\text{ -  border}$ 
\end{ptheorem}
\begin{ctheorem}
$ W(A)\supseteq conv(\sigma_\epsilon(A))\setminus \epsilon_1 \text{ -  border}$, here $\epsilon_1=\frac{2\epsilon}{1-\epsilon}\n{A}$.
\end{ctheorem}
\begin{proof} We first prove the following claim: 
$$conv(\sigma_\epsilon(a))\subseteq W(T)+B\left(0,\frac{2\epsilon}{1-\epsilon}\n{T}\right).
$$
If $\lambda \in \mathbb C$ and $d(\lambda,W(A))>0$, then $\lambda-A$ is invertible and (see \cite[Theorem 6.2-A]{taylor}, also Ex. 27.7, p501,\cite{lim}) $$\n{(\lambda I-A)^{-1}}\leq\dfrac{1}{d(\lambda,W(A))}.$$ Let $\lambda\in\sigma_\epsilon(A)$. If $\lambda \in \overline{W(A)}$, then the conclusion is obvious. Next, if $\lambda\in \sigma_\epsilon(A)\setminus \overline{W(A)}$, then 
\begin{align*}
d(\lambda,W(A))\leq\dfrac{1}{\n{(\lambda-A)^{-1}}}&\leq\epsilon \n{(\lambda-A)}
\leq\epsilon(\m{\lambda}+\n{A})\\
&\leq\epsilon\left(\dfrac{1+\epsilon}{1-\epsilon}\n{A}+\n{A}\right)
=\dfrac{2\epsilon}{1-\epsilon}\n{A}.
\end{align*}
Thus \[\sigma_{\epsilon}(A)\subseteq W(A)+B\left(0,\dfrac{2\epsilon}{1-\epsilon}\n{A}\right).\] Since the right hand side is a convex set we have established the claim. As $conv(\sigma_\epsilon(a))$ and $W(A)$ are convex sets we get the required result.
\end{proof}

Spectral mapping theorem describes the behavior of the spectrum under certain transformation. The 
following result belongs to that kind. It gives a precise information about the condition spectrum under linear transformation.
\begin{utheorem}
$\sigma(\alpha+\beta A)=\alpha+\beta\sigma(A)$ for all $\alpha, \beta \in \mathbb C$
\end{utheorem}
\begin{ptheorem}
$\Lambda_{\epsilon\m{\beta}}(\alpha+\beta A)=\alpha+\beta\Lambda_\epsilon(A)$ for all $\alpha, \beta \in \mathbb C$
\end{ptheorem}
\begin{ctheorem}
$\sigma_\epsilon(\alpha+\beta A)=\alpha+\beta\sigma_\epsilon(A)$ for all $\alpha, \beta \in \mathbb C$
\end{ctheorem}
\begin{proof}
For the case $\beta=0$ follows from the fact $\sigma_\epsilon(\alpha)=\{ \alpha \}$. Consider $\beta \neq 0$,
\begin{align*}
\n{z-(\alpha+\beta A)}\n{(z-(\alpha+\beta A))^{-1}}&=\n{\beta\dfrac{(z-\alpha)}{\beta}-\beta A}\n{\left(\beta\dfrac{(z-\alpha)}{\beta}-\beta A\right)^{-1}}\\
&=\n{\dfrac{(z-\alpha)}{\beta}- A}\n{\left(\dfrac{(z-\alpha)}{\beta}- A\right)^{-1}} \qedhere
\end{align*}
\end{proof}
The above results analyze the similarity between condition spectra and pseudospectra. The results in \cite{suku2, suku3} presents basic properties of conditionspectra and its connection with other areas. All these results together demonstrates the potential of condition spectrum and the need for further investigation. The computational aspect of this spectrum is yet to be investigated. There are other results in \cite{Mark} for which, at present, the generalizations in conditionspectrum are unknown.

%
%
%
\bibliography{newsukubib}

%
%
%
%
%

\end{document}